\documentstyle[12pt,amstex]{article}

\begin{document}
\newcommand{\bx}{\hfill\rule{.25cm}{.25cm}\medbreak}
\newtheorem{claim}{\indent Claim}
\newtheorem{theorem}{Theorem}[section]
\numberwithin{equation}{section}

\newtheorem{lemma}[theorem]{Lemma}
\newtheorem{proposition}[theorem]{Proposition}

\newtheorem{corollary}[theorem]{Corollary}
\newtheorem{definition}[theorem]{Definition}
\newtheorem{remark}[theorem]{Remark}
\newtheorem{remarks}[theorem]{Remarks}
\newtheorem{example}[theorem]{Example}

%\setcounter{section}{1}
%\setcounter{equation}{0}
%\setcounter{theorem}{0}

%\begin{theorem}\label{theorem1}
%\end{theorem}
%
%\begin{equation}
%\label{main-theorem-1}
%
%\end{equation}
%
%(\ref{main-theorem-1})
%
%By Theorem \ref{theorem1}.

\newtheorem{asum}{Assumption}

\newcommand{\bR}{\overline{R}}
\newcommand{\eg}{\frak g}
\newcommand{\stu}{{\frak stu}_n(R, -, \boldsymbol{\gamma})}
\newcommand{\eij}{E_{ij}(a)}
\newcommand{\ejk}{E_{jk}(b)}
\newcommand{\Tij}{T_{ij}(a, b)}
\newcommand{\eik}{E_{ik}(ab)}
\newcommand{\ekl}{E_{kl}(b)}
\newcommand{\st}  {{\frak st}_n(R)}
\newcommand{\stl}  {{\frak stl}_n(R)}
\newcommand{\sth}  {\widehat{{\frak stl}}_n(R)}
\newcommand{\stf} {{\frak stl}_4(R)}
\newcommand{\stfh} {\widehat{{\frak stl}}_4(R)}
\newcommand{\stft} {\widetilde{{\frak stl}}_4(R)}
\newcommand{\stfs} {{\frak stl}_4(R)^\sharp}
\newcommand{\stt} {{\frak stl}_3(R)}
\newcommand{\stth} {\widehat{{\frak stl}}_3(R)}
\newcommand{\sttt} {\widetilde{{\frak stl}}_3(R)}
\newcommand{\stts} {{\frak stl}_3(R)^\sharp}
\newcommand{\lij}{\gamma_i\gamma_j^{-1}}
\newcommand{\lji}{\gamma_j\gamma_i^{-1}}
\newcommand{\ljk}{\gamma_j\gamma_k^{-1}}
\newcommand{\gd}{\dot{\frak{g}}}
\def\ker{{\rm Ker\,}}
\def\Im{{\rm Im\,}}

\centerline{{\large \bf Second homology groups and universal
coverings of}}\centerline{{\large \bf Steinberg Leibniz algebras of
small characteristic}\,\footnote{This work is supported by NSF
grants 10571120, 10471096 of China and ``One Hundred Talents
Program" from University of Science and Technology of
China.\\\indent {\it 2000 Mathematics Subject Classification}:
17A32, 17B55, 17B60.}}

\centerline{(Appeared in {\it Comm.~Algebra} {\bf37} (2009), no. 2,
548-566.)}

 \vspace{1em} \centerline{\bf Qifen Jiang$^{\,1)}$, Ran Shen$^{\,1)}$ and Yucai
Su$^{\,2)}$
 }
\centerline{{\footnotesize $^{1)}$}\small\it Department of
Mathematics, Shanghai Jiaotong University, Shanghai 200240, China}
\centerline{\small\it E-mail: qfjiang@@sjtu.edu.cn}
\centerline{{\footnotesize $^{2)}$}\small\it Department of
Mathematics, University of Science and Technology of
China}\centerline{\small\it Hefei 230026,
China}\centerline{\small\it E-mail: ycsu@@ustc.edu.cn}

%{\bf(You need to check all Equation labels since I add some more
%equation numbers)}

%

{ \vspace{1.5em}
\begin{abstract}
It is known that the second Leibniz homology group $HL_2(\stl)$ of
the Steinberg Leibniz algebra $\stl$ is trivial for $n\geq 5$. In
this paper, we determine $HL_2(\stl)$ explicitly (which are shown to
be not necessarily trivial) for $n=3, 4$ without any assumption on
the base ring.
%As a consequence, we obtain $HL_2(sl_n(R))$ of
%matrix Leibniz algebra for $n=3, 4$.
\end{abstract}

\vspace{0.5em}

\noindent{\bf \S 1 Introduction}
\setcounter{section}{1}\setcounter{equation}{0}

 The concept of Leibniz algebras was introduced by Loday [Lo] in the study of Leibniz homology as a noncommutative analog
 of Lie algebras homology. A Leibniz algebra $L$ is a vector space equipped with a $K$-bilinear
 map [ , ]: $L\times L\longrightarrow L$ satisfying the Leibniz
 identity $[x,[y,z]]=[[x,y],z]-[[x,z],y]$ for all $x,y,z\in L$, where $K$ is a unital commutative ring.
 Clearly, a Lie algebra is a Leibniz algebra. For any Leibniz algebra $L$ there is an associated Lie algebra
 $L_{Lie}=L/\langle[x,x]\rangle$, where $\langle[x,x]\rangle$ is the two-sided ideal generated by all $[x,x], x\in L$.
 To study the second
 Leibniz homology group of Lie algebra  $sl_n(R)$ and Steinberg Lie algebra $\st$, Loday and Pirashvili [LP]
 introduced also the noncommutative Steinberg Leibniz algebra $\stl$, where $R$ is an associative algebra over
 a commutative ring $K$. In [L], Steinberg Leibniz algebra and its superalgebra were dicussed.
 Steinberg Lie algebras and their universal central extensions have been
 studied by many authors (e.g., [B1, KL, Ka, G1, G2, GS]). In most situations, the Steinberg Lie algebra $\st$ is
the universal central extension of the Lie algebra $sl_n(R)$ whose
kernel is isomorphic to the first cyclic homology group $HC_1(R)$ of
the associative algebra $R$ and the second Lie algebra homology
group $H_2(\st)=0$. In [Bl] and [KL], it was proved that
$H_2(\st)=0$ for $n\geq 5$. In [KL], it was mentioned without proof
that $H_2(\st)=0$ for $n=3, 4$ if $\frac{1}{2}\in$ the base ring
$K$. This was proved in [G1] for $n=3$ if $\frac{1}{6}\in K$ and for
$n=4$ if $\frac{1}{2}\in K$. Gao and Shang [GS] generalize the
result without any assumption on $K$. Since $\stl_{Lie}=\st$, it is
natural to consider the Second homology group $HL_2(\stl)$ of the
Steinberg Leibniz algebra $\stl$ in the category of Leibniz
algebras. In [LP], Loday and Pirashvili proved that $HL_2(\stl)=0$
for $n\geq 5$.

Motivated by [GS], in this paper, we will determine $HL_2(\stl)$
explicitly for $n=3, 4$ (which are not necessarily trivial). It is
equivalent to work on the Steinberg Leibniz algebras $\stl$ of
small characteristic for small $n$. This completes the
determination of the universal central extension of the Leibniz
algebras $\stl$ and $sl_n(R)$ as well. We would like to remark
that since skew-symmetry does not hold for Leibniz algebras, some
different approach to solve the problem seems to be necessary.
This is also one of our motivation to present this paper. The main
result in this paper is the following theorem (cf.~Theorem
$\ref{hl-1}$ and $\ref{hl-2}$), where the result for the first
case was obtained in [LP].
\begin{theorem} The second homology group of Steiberg Leibniz
algebra $\stl$ $($cf. Definition $\ref{main-def1})$ is
\begin{equation}
HL_2(\stl)=\cases 0 &\text{ if \ \ } n\geq 5,\\
R_2^6  &\text{ if \ \ } n=4,\\
R_3^6 &\text{ if \ \ } n=3,\endcases\notag
\end{equation}where $R_2^6, R_3^6$ are defined in Definition
$\ref{main-def2}$ and Definition $\ref{main-def3}.$
\end{theorem}

The paper is organized as follows. In Section 2, we review some
basic definitions and results on Steinberg Leibniz algebras
$\stl$. Section 3 will discuss the $n=4$ case. Section 4 will
handle the $n=3$ case.

\vspace{4mm}

\noindent{\bf \S 2 Preliminary }
\setcounter{section}{2}\setcounter{equation}{0}\setcounter{theorem}{0}
 Let $K$ be a
unital commutative ring.
\begin{definition} \rm A
{\it Leibniz algebra} $L$ is a vector space equipped with a
$K$-bilinear
 map $[\cdot  ,\cdot]:\, L\times L\rightarrow L,$  satisfying the Leibniz
 identity\begin{equation} \label{Leib-id}[x,[y,z]]=[[x,y],z]-[[x,z],y] \mbox{ \ \  for all \ } x,y,z\in L.
 \end{equation}
\end{definition}

 Clearly, a Lie algebra is a Leibniz algebra. For any Leibniz algebra $L$ there is an associated Lie algebra
 $L_{Lie}=L/\langle[x,x]\rangle$, where $\langle[x,x]\rangle$ is the two-sided ideal generated by all $[x,x], x\in L$.
\begin{definition}\rm Let $L$ be a Leibniz algebra over $K$, defined
the boundary map $\delta_n: L^{\otimes n}\rightarrow L^{\otimes
n-1}$ by $$\delta_n(x_{1}\otimes x_{2}\cdot\cdot\cdot\otimes
x_{n})=\sum\limits_{1\leq i<j\leq
n}(-1)^{j+1}x_{1}\otimes\cdot\cdot\cdot\otimes
x_{i-1}\otimes[x_{i},x_{j}]\otimes
x_{i+1}\otimes\cdot\cdot\cdot\otimes\hat{x_{j}}\otimes\cdot\cdot\cdot
x_{n},$$ where $\hat{x_{j}}$ means that the element $x_{j}$ is
omitted. The complex $( L^{\otimes n}, \delta )$ (where $L^{0}=K $
and $ \delta_{1}=0 )$ gives the Leibniz homology $HL_{\star}(L)$ of
the Leibniz algebra $L$ and $HL_n(L)=\ker\delta_{n}/
\Im\delta_{n+1}$ is called the $n$-th homology group of $L$.
 \end{definition}

Let $L$ be a Leibniz algebra over $K$, the {\it center} of $L$ is
defined to be $$\{ z\in L\,|\,[ z, L ]=[ L, z ]=0 \}.$$ A Leibniz
algebra L is called {\it perfect} if $[ L, L ]=L$. A {\it central
extension} of $L$ is a pair $(\hat{L}, \pi)$ where $\hat{L}$ is a
Leibniz algebra, and $\pi: \hat{L}\rightarrow L$ is a surjective
homomorphism such that $\ker\pi$ lies in the center of $\hat{L}$ and
the exact sequence
$$0\rightarrow\ker\pi\rightarrow\hat{L}\rightarrow L\rightarrow 0,$$
splits as $K$-modules. The pair $(\hat{L}, \pi)$ is a {\it universal
central extension} of $L$ if for every central extension
$(\tilde{L}, \tau)$ of $L$ there is a unique homomorphism $\psi:
\hat{L}\rightarrow \tilde{L}$ such that $\tau\circ\psi=\pi$.

The following result can be found in [LP].
\begin{proposition}The universal central extension of a Leibniz algebra $L$ exists if and only if $L$
is perfect. The universal central extension is unique up to
isomorphism. If $(\hat{L}, \pi)$ is the universal central extension
of $L$, then $HL_2(L)\cong\ker\pi$.
\end{proposition}

Let $R$ be a unital associative $K$-algebra. We always assume that
$R$ has a $K$-basis $\{r_\lambda\,|\,\lambda\in\Lambda\}$(where
$\Lambda$ is an index set), which contains the unit $1$ of $R$. We
denote by $gl_n(R)$ (where $n\geq 3$) the Leibniz algebra
consisting of all $n\times n$ matrices with coefficients in $R$
(which is in fact a Lie algebra), whose bracket is as
follows.$$[E_{ij}(a),
E_{kl}(b)]=\delta_{jk}E_{il}(ab)-\delta_{il}E_{kj}(ba),$$ for $a,
b\in R, 1\leq i, j, k, l\leq n$,  where $E_{ij}(x)$ is the matrix
with only non-zero element $x$ in position $(i, j)$.

The subalgebra $sl_n(R): =[gl_n(R), gl_n(R)]$ of $gl_n(R)$ is
generated by the elements $E_{ij}(a)$, $a\in R$, $1\leq i\neq
j\leq n$,  satisfying
\begin{eqnarray}&& [\eij , \ejk ] = \eik,
\\
&& [\eij , E_{ki}(b) ] =-E_{kj}(ba),
\end{eqnarray}
 for $i, j, k$  distinct and
 \begin{equation} [\eij , \ekl ] = 0 \mbox{ \ \ for \ }\ j\neq k,\ i \neq l.
 \end{equation}
\begin{definition}\mbox{\rm [LP]}\label{main-def1} \rm For $n\geq 3$, the {\it Steinberg Leibniz algebra} $\stl$
is a Leibniz algebra over $K$ defined by generators $X_{ij}(a)$,
$a\in R$, $1\leq i\neq j\leq n$, subject to the relations
\begin{eqnarray}
 &&X_{ij}(k_1a+k_2b)=k_1X_{ij}(a)+k_2X_{ij}(b)\ \  \mbox{  for }\  a, b\in R,\  k_1, k_2\in
 K,\\
&&[X_{ij}(a), X_{jk}(b)] = X_{ik}(ab),\\
&&[X_{ij}(a), X_{ki}(b)] = -X_{kj}(ba) \  \   \mbox{  for\ distinct
}\ \ i, j, k,\\
&&[X_{ij}(a), X_{kl}(b)] = 0\  \mbox{  for }\  j\neq k, i\neq l,
\end{eqnarray} where $ a, b\in R,$
$ 1\leq i,\  j,\  k,\  l \leq n$.
\end{definition}

One observes that both Leibniz algebras $sl_n(R)$ and $\stl$ are
perfect since $1\in R$. Define the Leibniz  algebra homomorphism
$\phi:\stl\rightarrow sl_n(R)$ by $\phi(X_{ij}(a))=E_{ij}(a)$.
Obviously, $\phi$ is an epimorphism.

The following result can be found in [LP, KL].
\begin{theorem} For $n\geq 3$ the kernel of $\phi$ is central in
$\stl$ and is isomorphic to $HH_1(R)$. Moreover if $n\geq 5$ then
$$0\rightarrow HH_1(R)\rightarrow\stl\rightarrow sl_n(R)\rightarrow 0$$
is the universal central extension of $sl_n(R)$ $($in the category
of Leibniz algebras$)$.
\end{theorem}

Here and below, $HH_1(R)$ denotes the {\it Hochschild homology
group} of $R$ with coefficients in $R$. This result means
$HL_2(\stl)=0$ for $n\geq 5$ and the universal central extension
of $sl_n(R)$ is also the universal central extension of $\stl$,
denoted by $\widehat{\frak{stl}}_n(R)$. Our purpose is to
determine $\widehat{\frak{stl}}_n(R)$ for any ring $K$ and $n\geq
3$. The following proposition can be similarly proved as in [AF]
for the Steinberg unitary Lie algebra case(Jacobi identity of Lie
algebra replaced by Leibniz identity).
\begin{lemma} Let $H:=\sum_{1\leq i\neq j\leq n}[X_{ij}(R),
X_{ji}(R)]$. Then $H$ is a subalgebra of $\stl$ containing the
center of $\stl$ such that $[H, X_{ij}(R)]\subseteq X_{ij}(R)$.
Moreover,
\begin{equation}\label{f4}\stl=H\oplus\sum_{1\leq i\neq j\leq
n}X_{ij}(R).\end{equation}
\end{lemma}

For a fixed $K$-basis $\{r_\lambda\,|\,\lambda\in\Lambda\}$ of $R$,
clearly
\begin{equation}\label{Gamma}\Gamma:=\{X_{ij}(r_\lambda)\,|\,\lambda\in\Lambda,\,1\leq i\neq
j\leq n\},\end{equation} is a $K$-basis of $\stl$.

%$\sum_{1\leq i\neq j\leq n}X_{ij}(R).$
%{\bf(note that I make a change here, please check if this change is
%right.)}

%
 We will see (cf.~Lemma
\ref{f5}) that the subalgebra $H$ has a more refined structure.
Setting
\begin{eqnarray}
&&T_{ij}(a, b) = [ X_{ij}(a), X_{ji}(b)],\\
&&t(a, b) = T_{1j}(a, b) - T_{1j}(ba, 1),
\end{eqnarray} for $a, b\in R, 1\leq i\neq j\leq n$. Note that $\Tij$
is $K$-bilinear, and so is $t(a,b)$. Further we have the following
\begin{lemma}\label{f3}For $a, b, c\in R$ and distinct $i, j, k$, we have
\begin{itemize}\parskip-3pt
 \item[\rm(1)]
$T_{ij}(a, bc)=T_{ik}(ab, c)+T_{kj}(ca,b)$,
 \item[\rm(2)] $T_{kj}(c,
1)+T_{jk}(c, 1)=0$,
 \item[\rm(3)] $t(a, b)$ \mbox{ does not depend on the choice of}
 $j$.
\end{itemize}\end{lemma} \noindent{\it Proof}.\ \ Using Leibniz
identity and the relations of generators of $\stl$, one has
\begin{eqnarray*}
T_{ij}(a, bc)&=&[X_{ij}(a), X_{ji}(bc)]
  =[X_{ij}(a), [X_{jk}(b),X_{ki}(c)]]\\
&=&[ [X_{ij}(a), X_{jk}(b)], X_{ki}(c)]-[[ X_{ij}(a), X_{ki}(c)],
X_{jk}(b)]\\
&=&[ X_{ik}(ab), X_{ki}(c)]-[-
X_{kj}(ca), X_{jk}(b)]\\
&=&T_{ik}(ab, c)+T_{kj}(ca, b).
\end{eqnarray*} This show that $(1)$
holds. By the $(1)$, we have \begin{equation}\label{f1}T_{ij}(ab,
c)=T_{ik}(a, bc)-T_{jk}(ca, b).\end{equation} Taking $b=1$ in
$(1)$ and $(\ref{f1})$, we have  $T_{ij}(a, c)=T_{ik}(a,
c)+T_{kj}(ca,1)$ and $T_{ij}(a, c)=T_{ik}(a, c)-T_{jk}(ca, 1)$.
Combining the two identities, we prove $(2)$.
 Taking $k\notin \{1, j\}$ and
using $(1)$, we have
\begin{eqnarray*}
t(a, b)&=&T_{1j}(a, b)-T_{1j}(ba, 1)\\
 &=&T_{1k}(a, b)+T_{kj}(ba, 1)-T_{1j}(ba, 1)\\
 &=&T_{1k}(a, b)-(T_{1j}(ba, 1)-T_{kj}(ba, 1))\\
 &=&T_{1k}(a, b)-T_{1k}(ba, 1),
 \end{eqnarray*}which proves $(3)$.
\hfill$\Box$\begin{lemma}\label{f5} Every $x\in H$ can be written as
the following form
$$ x = \sum_{i}t(a_i, b_i) + \sum_{2\leq j\leq n} T_{1j}(c_j, 1),$$
where $a_i, b_i, c_j\in R$.
\end{lemma}
\noindent{\it Proof}. ~~  Consider $T_{ij}(a, b)$ $($since $H$ is
generated by $T_{ij}(a, b), i\neq j\,)$. By Lemma $\ref{f3}$, we
have the following: Suppose $i=1$. Take $k\notin \{1, j\}$, then
\begin{eqnarray*}
T_{1j}(a, b)&=&T_{1k}(a, b)+T_{kj}(ba, 1)\\
 &=&T_{1k}(a, b)+T_{k1}(ba, 1)+T_{1j}(ba, 1)\\
 &=&T_{1k}(a, b)-T_{1k}(ba, 1)+T_{1j}(ba, 1))\\
 &=&t(a, b)+T_{1j}(ba, 1).
\end{eqnarray*}
Suppose $i, j\neq 1$. Then
\begin{eqnarray*}
T_{ij}(a, b)&=&T_{i1}(a b, 1)+T_{1j}(a, b)\\
 &=&-T_{1i}(a b, 1)+t(a, b)+T_{1j}(ba, 1)\\
 &=&t(a, b)-(T_{1i}(ab, 1)+T_{1j}(ba, 1)).
 \end{eqnarray*}
 Suppose $i\neq 1, j=1$. Take $k\notin \{1, i\}$, then
 \begin{eqnarray*}
T_{i1}(a, b)&=&T_{ik}(a , b)+T_{k1}(ba, 1)\\
 &=&T_{i1}(a b, 1)+T_{1k}(a, b)-T_{1k}(ba, 1)\\
 &=&t(a, b)-T_{1i}(ab, 1).
 \end{eqnarray*}This prove the lemma.\hfill$\Box$\vskip4pt

 For later use, we need the following results which are easy to
 check directly by Leibniz identity.
 \begin{lemma}\label{f6}For distinct $i, j, k, l$ and for $a, b, c\in
 R$, we have\begin{itemize}\parskip-1pt
 \item[\rm(1)] $[T_{ij}(a, b), X_{kl}(c)]=[X_{kl}(c), T_{ij}(a, b)]=0
 $,
 \item[\rm(2)] $[T_{ij}(a, b), X_{ik}(c)]=X_{ik}(abc)=-[X_{ik}(c), T_{ij}(a,
 b)]$,
 \item[\rm(3)] $[T_{ij}(a, b), X_{ki}(c)]=-X_{ki}(cab)=-[X_{ki}(c), T_{ij}(a,
 b)]$,
 \item[\rm(4)] $[T_{ij}(a, b), X_{kj}(c)]=X_{kj}(cba)=-[X_{kj}(c), T_{ij}(a,
 b)]$,
 \item[\rm(5)] $[T_{ij}(a, b), X_{ij}(c)]=X_{ij}(abc+cba)=-[X_{ij}(c), T_{ij}(a,
 b)]$,

 \item[\rm(6)] $[T_{ij}(a, b), X_{jk}(c)]=-X_{jk}(bac)=-[X_{jk}(c), T_{ij}(a,
 b)]$,
 \item[\rm(7)] $[t(a, b), X_{1i}(c)]=X_{1i}((ab-ba)c)=-[X_{1i}(c), t(a,
 b)]$,
 \item[\rm(8)] $[t(a, b), X_{i1}(c)]=-X_{i1}(c(ab-ba))=-[X_{i1}(c), t(a,
 b)]$,
 \item[\rm(9)] $[t(a, b), X_{jk}(c)]=0=-[X_{jk}(c), t(a, b)] $ for $j, k\geq
 2$.\end{itemize}
 \end{lemma}
\vspace{4mm}

\noindent{\bf \S 3 Universal central extension of $\stf$}
\setcounter{section}{3}\setcounter{equation}{0}\setcounter{theorem}{0}

In this section, we  determine the universal central extension
$\stfh$ of $\stf$ and compute $HL_2(\stf)$ without any assumption on
$R$.

For any positive integer $m$, let ${\mathcal I}_m$ be the ideal of
$R$ generated by the elements $ma$ and $ab-ba$, for $a,b \in R$.
One immediately has
\begin{lemma}\label{I-M}{\rm(}cf.~{\rm[GS])}
${\mathcal I}_m=mR+R[R,R] \text{ and } [R, R]R=[R, R]R$.
%\hfill$\Box$
\end{lemma}

Let $R_m:=R/{\mathcal I}_m$ be the quotient algebra over $K$ which
is commutative.  Write $\bar a=a+{\mathcal I}_m$ for $a\in R$.
Note that if $m=2$ then  $\overline{a}=-\overline{a}$ in $R_m$.
\begin{definition}\label{main-def2}
\rm We define ${\mathcal W}=R_2^6$ to be the direct sum of six
copies of $R_2$. For $1\leq m\leq 6$, we let
$\epsilon_{m}(\overline{a})=(0,\cdots,\overline{a},\cdots,0)$ be
the element of $\mathcal W$ such that the $m$-th coordinate is
$\overline{a}$ and zero otherwise.
\end{definition}

Let $S_4$ be the symmetric group of $\{1,2,3,4\}$. Let
$$P=\{(i,j,k,l)\,|\,\{i,j,k,l\}=\{1,2,3,4\}\},$$ be the set of all
the quadruples with the  distinct components. Then $S_4$ has a
natural transitive action on $P$ given by
$$\mbox{$\sigma((i,j,k,l))=(\sigma(i),\sigma(j),\sigma(k),\sigma(l))$ for
any $\sigma\in S_4$.}$$ Clearly $$G=\{(1),(13),(24),(13)(24)\},$$ is
a subgroup of $S_4$ with $[S_4:G]=6$. Then $S_4$ has a partition of
cosets with respect to $G$, denoted by
$$S_4=\bigsqcup_{m=1}^6\sigma_mG.$$ We  obtain a partition of $P$
by $$P=\bigsqcup_{m=1}^6P_m, \mbox{ \ where
}P_m=(\sigma_mG)((1,2,3,4)).$$ Define the index map
$\theta:P\rightarrow\{1,2,3,4,5,6\}$ by
$$\theta\left((i,j,k,l)\right)=m \text{ \ if } (i,j,k,l)\in P_m,\ \mbox{ \ for
}1\leq m\leq 6.$$ We fix $P_1=G((1,2,3,4))$, then we have
$(1,2,3,4)\in P_1$ and $\theta((1,2,3,4))=1$.

Using the decomposition $(\ref{f4})$ of $\stf$, we take $\Gamma$ as
in (\ref{Gamma}) with $n=4$.
%a $K$-basis $\Gamma$ of $\stf$, which
%contains $\{X_{ij}(r_\lambda)\,|\,\lambda\in\Lambda, 1\leq i\neq
%j\leq4\}$.
%{\bf(Note that I make some changes here)}
Define $\psi: \Gamma \times \Gamma \to \mathcal W$ by
$$\mbox{$\psi(X_{ij}(r),X_{kl}(s))=\epsilon_{\theta((i,j,k,l))}(\overline{rs})\in\mathcal
W$ for $r,s\in\{r_\lambda\,|\,\lambda\in\Lambda\}$ and distinct $i,
j, k, l$,}$$
 and $\psi=0$ otherwise. Then we obtain the $K$-bilinear
map $\psi:\stf\times\stf\to\mathcal W$ by linearity. We now have
\begin{lemma}\label{Le1} The bilinear map $\psi$ is a Leibniz $2$-cocycle.
\end{lemma}
\noindent{\it Proof}. It suffices to prove
\begin{equation}\label{J(x,y,z)}J(x,y,z):=
\psi(x,[y,z])+\psi([x,z],y)-\psi([x,y],z)=0,\end{equation} for any
$x,y,z\in\stf$.
%We need to check $J(x,y,z)=0$ for
%$x,y,z\in\Gamma$.
According to $(\ref{f4})$ and Lemma $\ref{f5}$,
the Steinberg Leibniz algebra $\stf$ has the decomposition :
\begin{equation} \stf= t(R,R)\oplus T_{12}(R,1)\oplus T_{13}(R,1)\oplus
T_{14}(R,1) \oplus \bigoplus_{1\leq i \neq j\leq n}X_{ij}(R),
\end{equation}
where $t(R,R)$ is the $K$-linear span of the elements $t(a,b)$.
Clearly, the number of elements of ${x,y,z}$ belonging to the
subalgebra $H$ such that $\psi([x,y],z)\neq 0$ is at most one. We
consider the following possibilities:

{\bf Case 1:}  Suppose there exists exactly one of $\{x,y,z\}$
belonging to $H$. Say, $x=X_{12}(a),y=X_{34}(b)$ and $z\in H$, where
$a,b\in R$ (we omit the other subcases since they are very similar,
although not identical). We can assume that either $z=t(c,d)$, where
$c,d\in R$, or $z=T_{1j}(c,1)$, where $2\leq j\leq 4$ and $c\in R$.

If $z=t(a,b)$, then according to the Leibniz identity and Lemma
$\ref{f6}$, we have
\begin{eqnarray}
J(x,y,z))&=&\psi([X_{12}(a),t(c,d)],X_{34}(b))\nonumber\\
&=&\psi(-X_{12}((cd-dc)a),X_{34}(b))\nonumber\\
&=&-\epsilon_1(\overline{(cd-dc)ab)}=0.\nonumber
\end{eqnarray}
If $z=T_{12}(c, 1)$, then
\begin{eqnarray}
J(x,y,z))&=&\psi([X_{12}(a),T_{12}(1,c)],X_{3,4}(b))\nonumber\\
&=&\psi(-X_{12}(ca+ac),X_{34}(b))\nonumber\\
&=-&\epsilon_1(\overline{(ca+ac)b})=0.\nonumber
\end{eqnarray}
If $z=T_{13}(c,1)$, then
\begin{eqnarray}
J(x,y,z))&=&\psi(X_{12}(a), [X_{34}(b),T_{13}(c, 1)])+\psi([X_{12}(a),T_{13}(c,1)],X_{34}(b))\nonumber\\
&=&\psi(X_{12}(a), X_{34}(cb))+\psi(-X_{12}(ca),X_{34}(b))\nonumber\\
&=&\epsilon_1(\overline{acb}-\overline{cab})=\epsilon_1(\overline{c(ab-ba)})=0.\nonumber
\end{eqnarray}
If $z=T_{14}(c,1)$, then
\begin{eqnarray}
J(x,y,z)&=&\psi(X_{12}(a),[X_{34}(b),T_{14}(c,1)])+\psi([X_{12}(a),T_{14}(c,1)],X_{34}(b))\nonumber\\
&=&\psi(X_{12}(a),-X_{34}(bc))+\psi(-X_{12}(ca),X_{34}(b))\nonumber\\
&=&-\epsilon_1(\overline{abc}+\overline{cab})=-\epsilon_1(\overline{abc-abc})=0.\nonumber
\end{eqnarray}

 {\bf Case 2:} Suppose none of
$\{x,y,z\}$ belongs to $H$. The nonzero terms of $J(x,y,z)$ must be
$\psi([X_{ik}(a),X_{kj}(b)],X_{kl}(c))$ or
$\psi([X_{il}(a),X_{lj}(b)],X_{kl}(c))$ for distinct $i,j,k,l$ and
$a,b,c\in R$.

In case $x=X_{ik}(a)$, $y=X_{kj}(b)$ and $z=X_{kl}(c)$, we have
\begin{eqnarray}
J(x,y,z))&=&\psi([X_{ik}(a),X_{kl}(c)],
X_{kj}(b))-\psi([X_{ik}(a),X_{kj}(b)], X_{kl}(c))\nonumber\\
&=&\psi(X_{il}(ac),X_{kj}(b))-\psi(X_{ij}(ab),X_{kl}(c))\nonumber\\
&=&\epsilon_{\theta((i,l,k,j))}(\overline{acb})-\epsilon_{\theta((i,j,k,l))}(\overline{abc})\nonumber\\
&=&\epsilon_{\theta((i,j,k,l))}(\overline{a(cb-bc)})=0.\nonumber
\end{eqnarray}

In case $x=X_{il}(a)$,$y=X_{lj}(b)$ and $z=X_{kl}(c)$, we have
\begin{eqnarray}
J(x,y,z))&=&\psi(X_{il}(a),[X_{lj}(b),
X_{kl}(c)])-\psi([X_{il}(a),X_{lj}(b)], X_{kl}(c))\nonumber\\
&=&\psi(X_{il}(a),-X_{kj}(cb))-\psi(X_{ij}(ab),X_{kl}(c))\nonumber\\
&=&-\epsilon_{\theta((i,l,k,j))}(\overline{acb})-\epsilon_{\theta((i,j,k,l))}(\overline{abc})\nonumber\\
&=&-\epsilon_{\theta((i,j,k,l))}(\overline{a(bc+cb)})=0,\nonumber
\end{eqnarray}
where the fourth equality follows from the fact that $(i,j,k,l)$ and
$(i,l,k,j)$ are in the same partition of $P$, i.e.
$\theta((i,j,k,l))=\theta((k,l,i,j))$. This is because that if
$(ijkl)=\sigma((1234))$ (where $\sigma\in S_4$), then
$(ilkj)=(\sigma\circ(24))((1234))$. The proof is completed.
\hfill$\Box$\vskip4pt

We therefore obtain a central extension of Leibniz algebra $\stf$:
\begin{equation}
0\rightarrow{\mathcal
W}\rightarrow\stfh\overset{\pi}\rightarrow\stf\rightarrow 0,
\end{equation}
i.e. \begin{equation} \stfh={\mathcal W}\oplus\stf,
\end{equation}with bracket
$$[(c,x),(c',y)]=(\psi(x,y),[x,y]),$$
for all $x,y\in\stf$ and $c,c'\in{\mathcal W}$, where $\pi:{\mathcal
W}\oplus\stf\rightarrow\stf$ is the second coordinate projection
map. Namely, $(\stfh,\pi)$ is a central extension of $\stf$. We will
show that $(\stfh,\pi)$ is the universal central extension of
$\stf$. To do this, we define a Leibniz algebra $\stfs$ to be the
Leibniz algebra generated by the symbols $X_{ij}^{\sharp}(a)$, $a\in
R$ and the $K$-linear space ${\mathcal W}$, satisfying the following
relations:
\begin{eqnarray}
\label{fen1}
&&X_{ij}{^\sharp}(k_1a+k_2b)=k_1X_{ij}{^\sharp}(a)+k_2X_{ij}{^\sharp}(b)
\mbox{ for }a, b\in R,\  k_1, k_2\in K,\\
\label{fen2}&&[X_{ij}^{\sharp}(a), X_{jk}^{\sharp}(b)]
=-[X_{jk}^{\sharp}(b), X_{ij}^{\sharp}(a)]= X_{ik}^{\sharp}(ab)
\text{ for distinct } i, j, k,\\
\label{fen3}&&[X_{ij}^{\sharp}(a),{\mathcal W}]=0=[{\mathcal W},
X_{ij}^{\sharp}(a)] \text{ for distinct } i, j,\\
\label{fen4}&&[X_{ij}^{\sharp}(a),X_{ij}^{\sharp}(b)]=0
\text{ for distinct } i, j,\\
\label{fen5}&&[X_{ij}^{\sharp}(a),X_{ik}^{\sharp}(b)]=0
\text{ for distinct } i, j, k,\\
\label{fen6}&&[X_{ij}^{\sharp}(a),X_{kj}^{\sharp}(b)]=0
\text{ for distinct } i, j, k,\\
\label{fen7}&&[X_{ij}^{\sharp}(a),
X_{kl}^{\sharp}(b)]=\epsilon_{\theta((i,j,k,l))}(\overline{ab})
\text{ for distinct } j, k, i, l,
\end{eqnarray}
where $a,b\in R$, $1\leq i,j,k, l\leq 4$. Since $1\in R$, we see
$\stfs$ is perfect. Clearly, there is a unique Leibniz algebra
homomorphism $\rho:\stfs\rightarrow\stfh$ such that
$\rho(X^\sharp_{ij}(a))=X_{ij}(a)$ and $\rho|_{\mathcal W}=id$. One
can easily observe that  $\rho$ is actually an isomorphism. Namely,

\begin{lemma}\label{qf7} $\rho:\stfs\rightarrow\stfh$ is a Leibniz algebra
isomorphism.
\end{lemma}
%{\bf Proof: }Let
%$T_{ij}^{\sharp}(a,b)=[X_{ij}^{\sharp}(a),X_{ji}^{\sharp}(b)]$.
%$Then one can easily check that for $a,b\in R$ and distinct
%$i,j,k$, one has
%\begin{equation}
%&T_{ij}^{\sharp}(a,bc)=T_{ik}^{\sharp}(ab,c)+T_{kj}^{\sharp}(ca,b)
%\end{equation}
%Indeed, the proof of (2.17) is the same as the proof in [KL] and
%[G1,Proposition 2.17].
%Put
%$t^{\sharp}(a,b)=T_{1j}^{\sharp}(a,b)-T_{1j}^{\sharp}(ba,1)$ for
%$a,b\in R, 2\leq j\leq 4$. Then $t^\sharp(a,b)$ does not depend on
%the choice of $j$. Also, one can easily check (as in [AF, Lemma
%1.1]) that
%$$\stfs=H^{\sharp}\oplus_{1\leq i \neq j\leq 4}
%X_{ij}^{\sharp}(R)$$ where
%$${\frak T}^{\sharp}= \left(\sum_{i, j, k, l \text{ are distinct}}[ X_{ij}^{\sharp}(R), X_{kl}^{\sharp}(R)]\right)
%\oplus\left(\sum_{1\leq i < j\leq 4}[ X_{ij}^{\sharp}(R),
%X_{ji}^{\sharp}(R)]\right).$$ It then follows from (2.16) and
%(2.17) above that
%\begin{equation}{\frak T}^{\sharp}={\mathcal W}\oplus\left(t^{\sharp}(R,R)\oplus T_{12}^{\sharp}(1,R)\oplus T_{13}^{\sharp}(1,R)\oplus
%T_{14}^{\sharp}(1,R)\right)\tag{2.18}
%\end{equation} where
%$t^{\sharp}(R,R)$ is the linear span of the elements
%$t^{\sharp}(a,b)$. So by Lemma 1.11, it suffices to show that the
%restriction of $\rho$ to $t^{\sharp}(R,R)$ is injective.

%Now the similar argument as  given in [AG,Lemma 6.18] shows that
 %there exists a linear map from $t(R,R)$ to $t^{\sharp}(R,R)$ so
%that $t(a,b)\mapsto t^{\sharp}(a,b)$ for $a,b\in R$. This map is
%the inverse of the restriction of $\rho$ to $t^{\sharp}(R,R)$.
%$\Box$

The analogue of the following theorem for Lie algebra was obtained
in [GS]. However, in our case, since the skew-symmetry does not hold
for Leibniz algebras, we need to find some different approach to
solve the problem.
\begin{theorem}\label{hl-1}$(\stfh,\pi)$ is the universal central extension of $\stf$ and
hence
$$HL_2(\stf)\cong\mathcal W.$$
\end{theorem}
{\it Proof. }
 Suppose
\begin{equation}0\rightarrow{\mathcal
V}\rightarrow\stft\overset{\tau}{\rightarrow}\stf\rightarrow\notag
0,\end{equation} is a central extension of $\stf$. We must show that
there exists a Leibniz algebra homomorphism
$\eta:\stfh\rightarrow\stft$ such that $\tau\circ\eta=\pi$. By Lemma
$\ref{qf7}$, it suffices to show that there exists a Leibniz algebra
homomorphism $\xi:\stfs\rightarrow\stft$ such that
$\tau\circ\xi=\pi\circ\rho$.

Using the $K$-basis $\{r_\lambda\ |\ \lambda\in\Lambda\}$ of $R$, we
choose a preimage $\widetilde{X}_{ij}(a)$ of $X_{ij}(a)$ under
$\tau$ for $1\leq i\neq j\leq 4$ and $a\in \{r_\lambda\ |\
\lambda\in\Lambda\}$. For distinct $i ,j, k$, let
$$[\widetilde{X}_{ik}(a),\widetilde{X}_{kj}(b)]=\widetilde{X}_{ij}(ab)+{\mu}_{ij}^k(a,b),$$
where ${\mu}_{ij}^k(a,b)\in{\mathcal V}$.  Similar to the
discussions in [GS] (replacing of the Jacobi identity by Leibniz
identity (\ref{Leib-id})), we obtain that ${\mu}_{ij}^k(a,b)$ is
independent of the choice of $k$ and by re-choosing the preimage
$\widetilde{X}_{ik}(a)$, we can suppose as in [GS],
\begin{equation}\label{f7}
[\widetilde{X}_{ik}(a),\widetilde{X}_{kj}(b)]=\widetilde{X}_{ij}(ab).
\end{equation}
We also need to consider the bracket
$[\widetilde{X}_{kj}(b),\widetilde{X}_{ik}(a)]$ for $a,b\in R$ and
distinct $i,j,k$ since there is no skew-symmetry for Leibniz
algebras. Our approach is different from that in [GS].
 Let
$$
[\widetilde{X}_{kj}(a),\widetilde{X}_{ik}(b)]=-\widetilde{X}_{ij}(ba)+\widetilde{\mu_{ij}}(a,b),$$
where $\widetilde{{\mu}_{ij}}(a,b)\in{\mathcal V}$ is independent of
$k$ as in ${\mu}_{ij}^k(a,b)$. Take distinct $ i,j,k,l$, then
\begin{equation}\label{UQ1}\left [[\widetilde{X}_{kl}(a),\widetilde{X}_{lj}(c)],\widetilde{X}_{ik}(b)\right
]=[\widetilde{X}_{kj}(ac),\widetilde{X}_{ik}(b)].\end{equation} The
left-hand side of (\ref{UQ1}) is, by Leibniz identity
(\ref{Leib-id}),
$$\left [\widetilde{X}_{kl}(a),[\widetilde{X}_{lj}(c),\widetilde{X}_{ik}(b)]]\right
]+\left
[[\widetilde{X}_{kl}(a),\widetilde{X}_{ik}(b)],\widetilde{X}_{lj}(c)\right
]=[-\widetilde{X}_{il}(ba),\widetilde{X}_{lj}(c)]=-\widetilde{X}_{ij}(bac),$$
since $[\widetilde{X}_{lj}(c),\widetilde{X}_{ik}(b)]\in \mathcal V$.
On the other hand,  the right-hand side of (\ref{UQ1}) is
$-\widetilde{X}_{ij}(bac)+\widetilde{{\mu}_{ij}}(ac,b)$. Thus
$\widetilde{{\mu}_{ij}}(ac,b)=0$. In particular, taking $c=1$, we
have $\widetilde{{\mu}_{ij}}(a,b)=0$. Therefore
\begin{equation}\label{f8}
[\widetilde{X}_{kj}(a),\widetilde{X}_{ik}(b)]=-\widetilde{X}_{ij}(ba),
\end{equation}for $a,b\in R$ and distinct $i,~j,~k$.
Now $(\ref{f7})$ and $(\ref{f8})$ imply $(\ref{fen2})$.

Next for $k\neq i , k\neq j$, we have
\begin{align}
[\widetilde{X}_{ij}(a),\widetilde{X}_{ij}(b)]&=\left[\widetilde{X}_{ij}(a),[\widetilde{X}_{ik}(b),\widetilde{X}_{kj}(1)]
\right]\notag\\
&=\left[[\widetilde{X}_{ij}(a),\widetilde{X}_{ik}(b)],\widetilde{X}_{kj}(1)\right]-
\left[[\widetilde{X}_{ij}(a),\widetilde{X}_{kj}(1)],\widetilde{X}_{ik}(b)\right]\notag\\
&=0-0=0,
\end{align}
as both of $[\widetilde{X}_{ij}(a),\widetilde{X}_{ik}(b)]$ and
$[\widetilde{X}_{ij}(a),\widetilde{X}_{kj}(1)]$ are in $\mathcal
V$. Thus, relation $(\ref{fen4})$ is obtained.

To obtain relation $(\ref{fen5})$, take $l\notin \{i,j,k\}$, then by
Leibniz identity (\ref{Leib-id}),
\begin{align}
[\widetilde{X}_{ij}(a),\widetilde{X}_{ik}(b)]&=\left[\widetilde{X}_{ij}(a),[\widetilde{X}_{il}(b),\widetilde{X}_{lk}(1)]\right]\notag\\
&=\left[[\widetilde{X}_{ij}(a),\widetilde{X}_{il}(b)],\widetilde{X}_{lk}(1)\right]-
\left[[\widetilde{X}_{ij}(a),\widetilde{X}_{lk}(1)],\widetilde{X}_{il}(b)\right]\notag\\
&=0-0=0,
\end{align}
since
$[\widetilde{X}_{ij}(a),\widetilde{X}_{il}(b)],[\widetilde{X}_{ij}(a),\widetilde{X}_{lk}(1)]\in
\mathcal V$. %Similarly, we have
%$$[\widetilde{X}_{ik}(b),\widetilde{X}_{ij}(a)]=0$$
Similarly, we have
\begin{equation}
[\widetilde{X}_{ij}(a),\widetilde{X}_{kj}(b)]=0,%=[\widetilde{X}_{kj}(b),\widetilde{X}_{ij}(a)]=0
\end{equation} for distinct $i,j,k$, i.e., we have relation  $(\ref{fen6})$.

To verify  $(\ref{fen7})$, set
$\widetilde{T}_{ij}(a,b)=[\widetilde{X}_{ij}(a),\widetilde{X}_{ji}(b)]$.
The following brackets are easily checked by the Leibniz identity.
\begin{align}\label{f9}&[\widetilde{T}_{ij}(a,b),\widetilde{X}_{ik}(c)]=\widetilde{X}_{ik}(abc)
=-[\widetilde{X}_{ik}(c),\widetilde{T}_{ij}(a,b)],\notag\\
&[\widetilde{T}_{ij}(a,b),\widetilde{X}_{kj}(c)]
=\widetilde{X}_{kj}(cba)=-[\widetilde{X}_{kj}(c),\widetilde{T}_{ij}(a,b)],\notag
 \text{and}\\
&[\widetilde{T}_{ij}(a,b),\widetilde{X}_{kl}(c)]=0=[\widetilde{X}_{kl}(c),\widetilde{T}_{ij}(a,b)].
\end{align}
Then we have
\begin{align}\label{f10}
[\widetilde{T}_{ij}(a,b),\widetilde{X}_{ij}(c)]&=\left[\widetilde{T}_{ij}(a,b),[\widetilde{X}_{ik}(c),\widetilde{X}_{kj}(1)]\right]\notag\\
&=\left[[\widetilde{T}_{ij}(a,b),\widetilde{X}_{ik}(c)],\widetilde{X}_{kj}(1)\right]-
\left[[\widetilde{T}_{ij}(a,b),\widetilde{X}_{kj}(1)],\widetilde{X}_{ik}(c)\right]\notag\\
&=\widetilde{X}_{ij}(abc)-(-\widetilde{X}_{ij}(cba))=\widetilde{X}_{ij}(abc+cba),
\end{align}for $a,b,c\in R$ and distinct $i,j,k,l$.
 Similarly $$[\widetilde{X}_{ij}(c),\widetilde{T}_{ij}(a,b)]=-\widetilde{X}_{ij}(abc+cba).$$

Next, for distinct $i,j,k,l$, let
$$[\widetilde{X}_{ij}(a),\widetilde{X}_{kl}(b)]=\nu^{ij}_{kl}(a,b),$$
where $\nu^{ij}_{kl}(a,b)\in \mathcal V$. By  $(\ref{f9})$ and
$(\ref{f10})$,
\begin{align}
2\nu^{ij}_{kl}(a,b)&=[\widetilde{X}_{ij}(2a),\widetilde{X}_{kl}(b)]
=\left[[\widetilde{T}_{ij}(1,1),\widetilde{X}_{ij}(a)],\widetilde{X}_{kl}(b)]\right]\notag\\
&=\left[\widetilde{T}_{ij}(1,1),[\widetilde{X}_{ij}(a),\widetilde{X}_{kl}(b)]\right]+
\left[[\widetilde{T}_{ij}(1,1),\widetilde{X}_{kl}(b)],\widetilde{X}_{ij}(a)]\right]\notag\\
&=0+0=0.\notag
\end{align}
So
\begin{equation}\label{f13}
\nu^{ij}_{kl}(a,b)=-\nu^{ij}_{kl}(a,b).
\end{equation}
Using Leibniz identity (\ref{Leib-id}), we have
\begin{eqnarray}
\nu^{il}_{kj}(bc,a)&=&\left
[[\widetilde{X}_{ik}(b),\widetilde{X}_{kl}(c)],\widetilde{X}_{kj}(a)\right]\nonumber\\
&=&\left[\widetilde{X}_{ik}(b),[\widetilde{X}_{kl}(c),\widetilde{X}_{kj}(a)]\right
]+\left
[[\widetilde{X}_{ik}(b),\widetilde{X}_{kj}(a)],\widetilde{X}_{kl}(c)\right]\nonumber\\
&=&[\widetilde{X}_{ij}(ba),\widetilde{X}_{kl}(c)]=\nu^{ij}_{kl}(ba,c).\nonumber
\end{eqnarray}
Taking $c=1$ and $b=1$ respectively, we have
$\nu^{il}_{kj}(b,a)=\nu^{ij}_{kl}(ba,1)$ and
$\nu^{il}_{kj}(c,a)=\nu^{ij}_{kl}(a,c)$. It then follows that
\begin{equation}\label{f12}\nu^{il}_{kj}(b,a)=\nu^{ij}_{kl}(a,b)=\nu^{ij}_{kl}(ba,1),
\end{equation}where $a,b\in R$ and $i,j,k,l$ are distinct.
So $\nu^{ij}_{kl}(R,R)=\nu^{ij}_{kl}(R,1)$.

Moreover, by $(\ref{f10})$, we get
\begin{align}\label{f11}
\nu^{ij}_{kl}(abc+cba,d)&=[\widetilde{X}_{ij}(abc+cba),\widetilde{X}_{kl}(d)]=\left
[[\widetilde{T}_{ij}(a,b),\widetilde{X}_{ij}(c)],\widetilde{X}_{kl}(d)\right
]\nonumber \\
&=\left[\widetilde{T}_{ij}(a,b),[\widetilde{X}_{ij}(c),\widetilde{X}_{kl}(d)]\right
]+\left
[[\widetilde{T}_{ij}(a,b),\widetilde{X}_{kl}(d)],\widetilde{X}_{ij}(c)\right]\nonumber \\
&=[\widetilde{T}_{ij}(a,b),\nu^{ij}_{kl}(c,d)]+[0,\widetilde{X}_{ij}(c)]=0.
\end{align}
Taking $c=d=1$ gives $\nu^{ij}_{kl}(ab+ba,1)=0$, i.e.
\begin{equation}\label{f15}\nu^{ij}_{kl}(a,b)=\nu^{ij}_{kl}(ba,1)=-\nu^{ij}_{kl}(ab,1)=\nu^{ij}_{kl}(ab,1)
=\nu^{ij}_{kl}(b,a).
\end{equation}
Letting $c=1$ in $(\ref{f11})$ and using $(\ref{f12})$ and
$(\ref{f13})$, we get
\begin{align}\label{f14}\nu^{ij}_{kl}(d(ab-ba),1)&=\nu^{ij}_{kl}(ab-ba,d)=\nu^{ij}_{kl}(ab,d)-\nu^{ij}_{kl}(ba,d)\notag\\
&=\nu^{ij}_{kl}(ab,d)+\nu^{ij}_{kl}(ba,d)\notag\\
&=\nu^{ij}_{kl}(ab+ba,d)=0,
\end{align}for $a,b,c,d\in R$ and distinct $i,j,k,l$.

 As for the other equalities of (\ref{fen7}), let
  $$[\widetilde{X}_{kl}(b),\widetilde{X}_{ij}(a)]=-\nu^{ij}_{kl}(a,b)+\nu',$$
where $\nu^{ij}_{kl}(a,b),\nu'\in \mathcal V$. By Leibniz identity
(\ref{Leib-id}) and $(\ref{f12})$, $(\ref{f13})$, we have
\begin{eqnarray}
\left[\widetilde{X}_{kl}(b),[\widetilde{X}_{ik}(a),\widetilde{X}_{kj}(1)]\right]&=&
\left[[\widetilde{X}_{kl}(b),\widetilde{X}_{ik}(a)],\widetilde{X}_{kj}(1)\right]-
\left[[\widetilde{X}_{kl}(b),\widetilde{X}_{kj}(1)],\widetilde{X}_{ik}(a)\right]\nonumber\\
&=&[-\widetilde{X}_{il}(ab),\widetilde{X}_{kj}(1)]-0\nonumber\\
&=&-\nu^{il}_{kj}(ab,1)=-\nu^{il}_{kj}(b,a)=-\nu^{ij}_{kl}(a,b).\nonumber
\end{eqnarray} So $\nu'=0$, and
\begin{equation}\label{f16}[\widetilde{X}_{kl}(b),\widetilde{X}_{ij}(a)]=-\nu^{ij}_{kl}(a,b).
\end{equation} Note that,
$(\ref{f13})$ and $(\ref{f14})$ show
\begin{equation}\nu^{ij}_{kl}({\mathcal I}_2,1)=0,\end{equation}
where ${\mathcal I}_2={\rm Span}\{2a,c(ab-ba)|a,b,c\in R)\}$
(cf.~Lemma $\ref{I-M}$). Moreover, for any $a\in R$, by
$(\ref{f13})$, $(\ref{f12})$, $(\ref{f15})$ and $(\ref{f16})$, we
get
$$\nu^{12}_{34}(a,1)=\nu^{14}_{32}(a,1)=\nu^{34}_{12}(a,1)=\nu^{32}_{14}(a,1).$$
It shows that the subgroup $G=\{(1),(13),(24),(13)(24)\}$ of $S_4$
fixes $\nu^{12}_{34}(a,1)$. Now similar arguments as in [GS]
complete the proof of the theorem.  \hfill$\Box$

\vspace{4mm}

\noindent{\bf \S 4 The second homology group $HL_2(\stt)$ of $\stt$}
\setcounter{section}{4}\setcounter{equation}{0}\setcounter{theorem}{0}

In this section we compute $HL_2(\stt)$. Recall from Lemma \ref{I-M}
that ${\mathcal I}_3=3R+R[R,R]$, and $R_3=R/{\mathcal I}_3$ is an
associative commutative $K$-algebra.

\begin{definition}\label{main-def3}Denote
${\mathcal U}=R_3^6$, and we also use $R_3^{(i)},i=-3,-2,\cdots,3$
to denote a copy of $R_3$. For $\overline{a}\in R_3$,
$\overline{a}^{(i)}$ will denote the corresponding element
$(0,\cdots,\overline{a},\cdots,0)$ in ${\mathcal U}$.
\end{definition}

For convenience, for $1\leq m\neq n\leq 3$ we use the symbol:
$${\rm sign}(m,n)=\cases
1,&\text { if } m<n,\\-1,&\text { if } m>n.\endcases $$
% Taking the $K$-basis $\Gamma$ of $\stt$:
%$\Gamma=\{X_{ij}(r_\lambda)\ |\ \lambda\in\Lambda,1\leq i\neq
%j\leq 3\}$.
Take $\Gamma$ as in (\ref{Gamma}) with $n=3$.
%
%{\bf(note that I make some change here)}
%
We define $\psi:\Gamma\times\Gamma\rightarrow\mathcal U$ by
$$\psi(X_{ij}(r),X_{ik}(s))={\rm sign}(j,k)(\overline{rs})^{(i)}$$
$$\psi(X_{ij}(r),X_{kj}(s))={\rm sign}(i,k)(\overline{rs})^{(-j)},$$ for $r,s\in
\{r_\lambda\ |\ \lambda\in\Lambda\}$ and distinct $i,j,k$, and
$\psi(x,y)=0 \text { otherwise. }$ Then $\psi$ can be extended to
a bilinear map $\stt\times\stt\rightarrow\mathcal U$. We have
\begin{lemma} The bilinear map $\psi$ is a Leibniz $2$-cocycle.
\end{lemma}
{\bf Proof: } Similar to the proof of Lemma $\ref{Le1}$, we show
$J(x,y,z)=0$ for $x,y,z\in\stt$ (cf.~(\ref{J(x,y,z)})). We have
the decomposition,
\begin{equation}\label{dec-id}
\stt=t(R,R)\oplus T_{12}(R,1)\oplus T_{13}(R,1)\oplus
\bigoplus_{1\leq i \neq j\leq 3}X_{ij}(R).\end{equation} As in the
proof of Lemma $\ref{Le1}$, we can suppose at most one of
$\{x,y,z\}$ is in the subalgebra $H$.
%{\bf(note that I change G to $H$. Same below.)}
%G.
%
We consider the following cases.

{\bf Case 1:} Suppose $z\in H$. %$z\in G$.

We first verify two subcases $x=X_{12}(a),y=X_{13}(b)$ or
$x=X_{21}(a),y=X_{31}(b)$ for $a,b\in R$.  By (\ref{dec-id}), we
may assume that either $z=t(c,d)$, where $c,d\in R$, or
$z=T_{1j}(c,1)$, where $2\leq j\leq 3$ and $c\in R$.

In case $x=X_{12}(a),y=X_{13}(b)$, we have, according to the
Leibniz identity, if $z=t(c,d)$, then
\begin{eqnarray}
J(x,y,z))&=&\psi(X_{12}(a),[X_{13}(b),t(c,d)])+\psi([X_{12}(a),t(c,d)],X_{13}(b))\nonumber\\
&=&-\psi(X_{12}(a),X_{13}((cd-dc)b))-\psi(X_{12}((cd-dc)a),X_{13}(b))\nonumber\\
&=&-(\overline{a(cd-dc)b+(cd-dc)ab})^{(1)}=0.\nonumber
\end{eqnarray}
If $z=T_{12}(c,1)$, then
\begin{eqnarray}
J(x,y,z))&=&\psi(X_{12}(a),[X_{13}(b),T_{12}(c,1)])+\psi([X_{12}(a),T_{12}(c,1)],X_{13}(b))\nonumber\\
&=&\psi(X_{12}(a), -X_{13}(cb))+\psi(-X_{12}(ca+ac),X_{13}(b))\nonumber\\
&=&-(\overline{acb+(ca+ac)b})^1=-(\overline{3abc})^{(1)}=0.\nonumber
\end{eqnarray}
If $z=T_{13}(c,1)$, then
\begin{eqnarray}
J(x,y,z))&=&\psi(X_{12}(a),[X_{13}(b),T_{13}(c,1)])+\psi([X_{12}(a),T_{13}(c,1)],X_{13}(b))\nonumber\\
&=&\psi(X_{12}(a), -X_{13}(cb+bc))+\psi(-X_{12}(ca),X_{13}(b))\nonumber\\
&=&-(\overline{a(cb+bc)+cab})^1=-(\overline{3abc})^{(1)}=0.\nonumber
\end{eqnarray}
In case $x=X_{21}(a),y=X_{31}(b)$, if $z=t(c,d)$, then
\begin{eqnarray}
J(x,y,z))&=&\psi(X_{21}(a),[X_{31}(b),t(c,d)])+\psi([X_{21}(a),t(c,d)],X_{31}(b))\nonumber\\
&=&\psi(X_{21}(a), X_{31}(b(cd-dc)))+\psi(X_{21}(a(cd-dc)),X_{31}(b))\nonumber\\
&=&(\overline{(ab(cd-dc)+a(cd-dc)b)})^{(-1)}=0.\nonumber
\end{eqnarray}
If $z=T_{12}(c,1)$, then
\begin{eqnarray}
J(x,y,z))&=&\psi(X_{21}(a),[X_{31}(b),T_{12}(1,c)])+\psi([X_{21}(a),T_{12}(c,1)],X_{31}(b))\nonumber\\
&=&\psi(X_{21}(a),X_{31}(bc))+\psi(X_{21}(ca+ac),X_{31}(b))\nonumber\\
&=&(\overline{abc+(ca+ac)b})^{(-1)}=(\overline{3abc})^{(-1)}=0.\nonumber
\end{eqnarray}
If $z=T_{13}(c,1)$, then
\begin{eqnarray}
J(x,y,z))&=&\psi(X_{21}(a),[X_{31}(b),T_{13}(c,1)])+\psi([X_{21}(a),T_{13}(1,c)],X_{31}(b))\nonumber\\
&=&\psi(X_{21}(a),X_{31}(cb+bc))+\psi(X_{21}(ac),X_{31}(b))\nonumber\\
&=&(\overline{a(cb+bc)+acb})^{(-1)}=(\overline{3abc})^{(-1)}=0.\nonumber
\end{eqnarray}

As for the other subcases, they are similar to the above subcases
except the following subcase: $x=X_{23}(a), y=X_{21}(b)$ and $z=t(c,
d)$. In this situation, we have
\begin{eqnarray}
J(x,y,z))&=&\psi(X_{23}(a),[X_{21}(b),t(c,d)])+\psi([X_{23}(a),t(c,d)],X_{21}(b))\nonumber\\
&=&\psi(X_{23}(a),X_{21}(b(cd-dc)))+\psi(0,X_{21}(b))\nonumber\\
&=&-(\overline{ab(cd-dc)})^{(2)}=0.\nonumber
\end{eqnarray}

 {\bf Case 2:} Suppose
there is none of $\{x,y,z\}$ belonging to $H$. %$G$.
Then the
nontrivial terms of $J(x,y,z)$ must be
$\psi([X_{ik}(a),X_{kj}(b)],X_{ik}(c))$ or
$\psi([X_{ik}(a),X_{kj}(b)],X_{kj}(c))$ for $a,b,c\in R$ and
distinct $i,j,k$.

If: $x=X_{ik}(a)$, $y=X_{kj}(b)$ and $z=X_{ik}(c)$, then
\begin{eqnarray}
J(x,y,z))&=&\psi(X_{ik}(a),[X_{kj}(b),X_{ik}(c)])-\psi([X_{ik}(a),X_{kj}(b)],X_{ik}(c))\nonumber\\
&=&\psi(X_{ik}(a),-X_{ij}(cb))-\psi(X_{ij}(ab),X_{ik}(c))\nonumber\\
&=&-{\rm sign}(k,j)(\overline{acb-abc})^{(i)}=0.\nonumber
\end{eqnarray}
If $x=X_{ik}(a)$, $y=X_{kj}(b)$ and $z=X_{kj}(c)$, then
\begin{eqnarray}
J(x,y,z))&=&\psi([X_{ik}(a),X_{kj}(c)],X_{kj}(b))-\psi([X_{ik}(a),X_{kj}(b)],X_{kj}(c))\nonumber\\
&=&\psi(X_{ij}(ac),X_{kj}(b))-\psi(X_{ij}(ab),X_{kj}(c))\nonumber\\
&=&{\rm sign}(i,k)(\overline{acb-abc})^{(-j)}=0.\nonumber
\end{eqnarray}
The proof is completed.    \hfill$\Box$\vskip4pt

It is similar to the $\stf$ case, we obtain a central extension of
$\stt$:
\begin{equation}
0\rightarrow{\mathcal
U}\rightarrow\stth\overset{\pi}\rightarrow\stt\rightarrow 0,
\end{equation}
i.e. \begin{equation} \stth={\mathcal U}\oplus\stt,
\end{equation}
and define $\stts$ to be the Leibniz algebra generated by the
symbols $X_{ij}^{\sharp}(a)$, $a\in R$ and the $K$-linear space
${\mathcal U}$, satisfying the following relations:
\begin{align}\label{n=3-1} &X_{ij}{^\sharp}(k_1a+k_2b)=k_1X_{ij}{^\sharp}(a)+k_2X_{ij}{^\sharp}(b)
\mbox{ for }a, b\in R,\  k_1, k_2\in
K,\\
\label{n=3-2}
&[X_{ij}^{\sharp}(a), X_{jk}^{\sharp}(b)]=-[X_{jk}^{\sharp}(b),X_{ij}^{\sharp}(a)] = X_{ik}^{\sharp}(ab) \text{ for distinct } i, j, k, \\
\label{n=3-3}
&[X_{ij}^{\sharp}(a),{\mathcal U}]=0=[{\mathcal U},X_{ij}^{\sharp}(a)] \text{ for distinct } i, j, \\
\label{n=3-4}
&[X_{ij}^{\sharp}(a),X_{ij}^{\sharp}(b)]=0 \text{ for distinct } i, j, \\
\label{n=3-5}&[X_{ij}^{\sharp}(a),X_{ik}^{\sharp}(b)]={\rm sign}(j,k)(\overline{ab})^i%=-[X_{ik}^{\sharp}(b)(a),X_{ij}^{\sharp}],
\text{ for distinct } i, j, k, \\
\label{n=3-6} &[X_{ij}^{\sharp}(a),X_{kj}^{\sharp}(b)]={\rm
sign}(i,k)(\overline{ab})_j%=-[X_{kj}^{\sharp}(b),X_{ij}^{\sharp}(a)],
\text{ for distinct } i, j, k,
\end{align}
where $a,b\in R$, $1\leq i,j,k\leq 3$ are distinct. Then $\stts$
is perfect and there is a unique Leibniz algebra homomorphism
$\rho:\stts\rightarrow\stth$. %We have separated the case (1.5)
%into three subcases (3.9)-(3.11).

As in Lemma $\ref{qf7}$, we have
\begin{lemma} $\rho:\stts\rightarrow\stth$ is a Leibiz algebra isomorphism. %$\Box$
\end{lemma}

Now we can state  the main theorem of this section.
\begin{theorem}\label{hl-2}$(\stth,\pi)$ is the universal central extension of $\stt$ and
hence
$$HL_2(\stt)\cong\mathcal U.$$
\end{theorem}
{\it Proof. } The idea to prove this theorem is similar to that in
the proof of Theorem 3.5. But there are some slight differences. The
point is that since $1\leq i, j, k\leq 3$, if $i, j, k$ are
distinct, then $k$ is uniquely determined once $i, j$ are chosen.

Suppose
\begin{equation}0\rightarrow{\mathcal
V}\rightarrow
\sttt\overset{\tau}{\rightarrow}\stt\rightarrow\notag
0,\end{equation} is a central extension of $\stt$. We must show
that there exists a Leibniz algebra homomorphism
$\eta:\stth\rightarrow\sttt$ so that $\tau\circ\eta=\pi$. Thus, by
Lemma 4.3, it suffices to show that there exists a Leibniz algebra
homomorphism $\xi:\stts\rightarrow\sttt$ such that
$\tau\circ\xi=\pi\circ\rho$.
 Choose a preimage $\widetilde{X}_{ij}(a)$ of $X_{ij}(a)$ as in
 Section 3,
we need to check that relations (\ref{n=3-1})--(\ref{n=3-6}) are
satisfied.

Again set
$$\begin{array}{ll}
\widetilde{T}_{ij}(a,b)=[\widetilde{X}_{ij}(a),\widetilde{X}_{ji}(b)],
\\[4pt]{}
[\widetilde{X}_{ik}(a),\widetilde{X}_{kj}(b)]=\widetilde{X}_{ij}(ab)+\mu_{ij}(a,b),
\\[4pt]{}
[\widetilde{X}_{kj}(b),\widetilde{X}_{ik}(a)]=-\widetilde{X}_{ij}(ab)+\widetilde{\mu_{ji}}(b,a),\end{array}$$
for $a,b\in R$ and distinct $i, j, k$, where
${\mu}_{ij}(a,b),\widetilde{\mu_{ji}}(b,a)\in{\mathcal V}$. Then by
Leibiz identity and the definition of $\widetilde{T}_{ij}(a,b)$, we
have
$$[\widetilde{T}_{ij}(a,b),\widetilde{X}_{ik}(c)]=\widetilde{X}_{ik}(abc)+{\mu}_{ik}(a, bc).$$
On the other hand, we also have
\begin{eqnarray*}
[\widetilde{T}_{ij}(a,b),\widetilde{X}_{ik}(c)]&=&\left[\widetilde{T}_{ij}(a,b),[\widetilde{X}_{ij}(1),\widetilde{X}_{jk}(c)]
\right]\\&=&\left[[\widetilde{T}_{ij}(a,b),\widetilde{X}_{ij}(1)],\widetilde{X}_{jk}(c)\right]
-\left[[\widetilde{T}_{ij}(a,b),\widetilde{X}_{jk}(c)],\widetilde{X}_{ij}(1)
\right]\\&=&[\widetilde{X}_{ij}(ab+ba),\widetilde{X}_{jk}(c)]-[-\widetilde{X}_{jk}(bac),\widetilde{X}_{ij}(1)]
\\&=&\widetilde{X}_{ik}(abc)+\widetilde{X}_{ik}(bac)+{\mu}_{ik}(ab,
c)+{\mu}_{ik}(ba,
c)-\widetilde{X}_{ik}(bac)+\widetilde{\mu_{ki}}(bac,1)
\\&=&\widetilde{X}_{ik}(abc)+{\mu}_{ik}(ab, c)+{\mu}_{ik}(ba,
c)+\widetilde{\mu_{ki}}(bac,1).
\end{eqnarray*} So we get
\begin{equation}\label{qf8}{\mu}_{ik}(a, bc)={\mu}_{ik}(ab, c)+{\mu}_{ik}(ba,
c)+\widetilde{\mu_{ki}}(bac,1).
\end{equation} Taking $b=1$ in $(\ref{qf8})$ gives
\begin{equation}{\mu}_{ik}(a,c)+\widetilde{\mu_{ki}}(ac,1)=0.\end{equation}
It follows that
\begin{equation}\label{qf9}{\mu}_{ik}(a,c)=-\widetilde{\mu_{ki}}(ac,1)={\mu}_{ik}(1,ac)={\mu}_{ik}(ac,1).\end{equation}
Similarly, by Leibniz  identity, we have
\begin{eqnarray*}
\left[\widetilde{T}_{ij}(1,1),[\widetilde{X}_{ij}(a),\widetilde{X}_{jk}(b)]\right]
&=&[\widetilde{T}_{ij}(1,1),\widetilde{X}_{ik}(ab)]=\widetilde{X}_{ik}(ab)+{\mu}_{ik}(1,ab)\\
&=&\left[[\widetilde{T}_{ij}(1,1),\widetilde{X}_{ij}(a)],\widetilde{X}_{jk}(b)\right]
-\left[[\widetilde{T}_{ij}(1,1),\widetilde{X}_{jk}(b)],\widetilde{X}_{ij}(a)\right]\\
&=&\widetilde{X}_{ik}(ab)+2{\mu}_{ik}(a,b)+\widetilde{\mu_{ki}}(b,a),
\end{eqnarray*}
which yields\begin{equation}\label{qf10}
{\mu}_{ik}(1,ab)=2{\mu}_{ik}(a,b)+\widetilde{\mu_{ki}}(b,a).\end{equation}
From $(\ref{qf9})$ and $(\ref{qf10})$, we have
\begin{equation}{\mu}_{ik}(a,b)+\widetilde{\mu_{ki}}(b,a)=0.\end{equation}
Now replacing $\widetilde{X}_{ij}(a)$ by
$\widetilde{X}_{ij}(a)+\mu_{ij}(1,a)$ which satisfies (\ref{n=3-1}),
we have at once
\begin{align}
[\widetilde{X}_{ik}(a),\widetilde{X}_{kj}(b)]=\widetilde{X}_{ij}(ab),\\
[\widetilde{X}_{kj}(b),\widetilde{X}_{ik}(a)]=-\widetilde{X}_{ij}(ba),
\end{align} for $a,b\in R$, distinct $i, j, k$.
This gives (\ref{n=3-2}).  The proof of  relation (\ref{n=3-4}) is
exactly the same as (\ref{qf10}).

To show $\widetilde{X}_{ij}(a)$ satisfies (\ref{n=3-5}) and
(\ref{n=3-6}), we define
$$[\widetilde{X}_{ij}(a),\widetilde{X}_{ik}(b)]=\nu^i_{jk}(a,b) \ \text{ and }\  [\widetilde{X}_{ij}(a),\widetilde{X}_{kj}(b)]
=\nu^{ik}_j(a,b),$$ where $\nu^i_{jk}(a,b),\nu^{ik}_j(a,b)\in
\mathcal V$. Then
\begin{align}
\nu^i_{jk}(a,b)&=[\widetilde{X}_{ij}(a),\widetilde{X}_{ik}(b)]=
\left[\widetilde{X}_{ij}(a),[\widetilde{X}_{ij}(1),\widetilde{X}_{jk}(b)]\right]\notag\\
&=-[\widetilde{X}_{ik}(ab),\widetilde{X}_{ij}(1)]=-\nu^i_{kj}(ab,1).
\end{align} and
\begin{align}
\nu^i_{jk}(a,b)&=[\widetilde{X}_{ij}(a),\widetilde{X}_{ik}(b)]=
\left[[\widetilde{X}_{ik}(a),\widetilde{X}_{kj}(1)],\widetilde{X}_{ik}(b)\right]\notag\\
&=[\widetilde{X}_{ik}(a),-\widetilde{X}_{ij}(b)]=-\nu^i_{kj}(a,b)
\notag
\\
&=\left[[\widetilde{X}_{ik}(1),\widetilde{X}_{kj}(a)],\widetilde{X}_{ik}(b)\right]\notag\\
&=[\widetilde{X}_{ik}(1),-\widetilde{X}_{ij}(ba)]=-\nu^i_{kj}(1,ba).
\end{align} So we have
\begin{equation}\label{qf1} \nu^i_{jk}(a,b)=-\nu^i_{kj}(a,b)=-\nu^i_{kj}(ab,1)=-\nu^i_{kj}(1,ba).\end{equation}
Moreover, we have
\begin{align}
0&=\left[\widetilde{T}_{ij}(1,1),[\widetilde{X}_{ij}(a),\widetilde{X}_{ik}(1)]\right]\notag\\
&=\left[[\widetilde{T}_{ij}(1,1),\widetilde{X}_{ij}(a)],\widetilde{X}_{ik}(1)]\right]-
\left[[\widetilde{T}_{ij}(1,1),\widetilde{X}_{ik}(1)],\widetilde{X}_{ij}(a)\right]\notag\\
&=2[\widetilde{X}_{ij}(a),\widetilde{X}_{ik}(1)]-[\widetilde{X}_{ik}(1),\widetilde{X}_{ij}(a)]\notag\\
&=2\nu^i_{jk}(a,1)-\nu^i_{kj}(1,a)\notag\\
&=2\nu^i_{jk}(1,a)+\nu^i_{jk}(1,a)=\nu^i_{jk}(1,3a).\label{qf2}
\end{align}
Similarly, we have
\begin{equation}\label{qf3}
\nu^{ik}_j(a,b)=\nu^{ik}_j(b,a)=\nu^{ik}_j(ab,1)=-\nu^{ki}_j(a,b)=-\nu^{ki}_j(ba,1),
\end{equation} and
\begin{equation}\label{qf4} \nu^{ik}_j(1,3a)=0. \end{equation}
Set
$\widetilde{t}(a,b)=\widetilde{T}_{1j}(a,b)-\widetilde{T}_{1j}(ba,1)$,
which does not depend on the choice of $j$. By the Leibniz
identity and $(\ref{qf1})$--$(\ref{qf4})$, we have
\begin{align}
0&=\left[\widetilde{t}(a,b),[\widetilde{X}_{12}(1),\widetilde{X}_{13}(c)]\right]\notag\\
&=\left[[\widetilde{t}(a,b),\widetilde{X}_{12}(1)],\widetilde{X}_{13}(c)]\right]-
\left[[\widetilde{t}(a,b),\widetilde{X}_{13}(c)],\widetilde{X}_{12}(1)\right]\notag\\
&=[\widetilde{X}_{12}(ab-ba),\widetilde{X}_{13}(c)]-[\widetilde{X}_{13}((ab-ba)c),\widetilde{X}_{12}(1)]\notag\\
&=\nu^1_{23}(ab-ba,c)-\nu^1_{32}((ab-ba)c,1)=\nu^1_{23}((ab-ba)c,1)+\nu^1_{23}((ab-ba)c,1)\notag\\
&=\nu^1_{23}(2(ab-ba)c,1)=\nu^1_{23}(1,2(ab-ba)c)=-\nu^1_{23}(1,(ab-ba)c),\notag
\end{align}
i.e. $\nu^1_{23}(1,(ab-ba)c)=0$. We also have
\begin{align}
0&=\left[\widetilde{t}(a,b),[\widetilde{X}_{12}(1),\widetilde{X}_{32}(c)]\right]\notag\\
&=\left[[\widetilde{t}(a,b),\widetilde{X}_{12}(1)],\widetilde{X}_{32}(c)]\right]\notag\\
&=[\widetilde{X}_{12}(ab-ba),\widetilde{X}_{32}(c)]\notag\\
&=\nu^{13}_{2}(ab-ba,
c)=\nu^{13}_{2}(c,ab-ba)=\nu^{13}_{2}((ab-ba)c,1)=\nu^{13}_{2}(1,(ab-ba)c).\notag
\end{align}
More generally,
\begin{equation}\nu^i_{jk}(1,(ab-ba)c)=0 \text {     and    }
\nu^{ik}_j(1,(ab-ba)c)=0,
\end{equation}
for $a,b,c\in R$ and distinct $1\leq i,j,k\leq 3$. Above
discussions prove
\begin{equation}\nu^i_{jk}(1, {\mathcal I}_3)=0 \text { and }
\nu^{ik}_j(1, {\mathcal I}_3)=0.
\end{equation} Now similar arguments as in [GS] complete the proof of the
theorem.\hfill$\Box$

\end{document}